\documentclass[11pt, a4paper]{article}

\usepackage[utf8]{inputenc}
\usepackage{alphalph, amsfonts, amsmath, amssymb, mathtools, siunitx}
\usepackage[margin=0.85in]{geometry}

\newtheorem{thrm}{Theorem}[section]
\newtheorem{lmm}[thrm]{Lemma}

\newtheorem{proposition}[thrm]{Proposition}

\newtheorem{assumption}[thrm]{Assumption}

\newenvironment{proof}{\noindent {\bf Proof. }}{\hfill $\Box$ \newline\par}

\newcommand{\Pmax}{P_{\max}}

\newcommand{\wmax}{w^{\max}}
\newcommand{\wmin}{w^{\min}}
\newcommand{\winit}{w_{\text{init}}}
\newcommand{\wfin}{w_{\text{fin}}}

\usepackage{algorithm}
\usepackage{algpseudocode}

\title{Hidden convexity property of a speed planning problem}
\author{Stefano Ardizzoni, Luca Consolini, Mattia Laurini, Marco Locatelli
\thanks{All authors are with the University of Parma, Department of Engineering and Architecture, Parco Area delle Scienze 181/A, 43124 Parma, Italy. E-mails:
{\tt\footnotesize \{stefano.ardizzoni, luca.consolini, mattia.laurini, marco.locatelli\}@unipr.it}
}}

\date{}
\begin{document}
\maketitle
\begin{abstract}
In this paper we address the speed planning problem for a vehicle along a predefined path.
A weighted average of two (conflicting) terms, energy consumption and travel time, is minimized.
After deriving a non-convex mathematical model of the problem, we introduce a convex relaxation of the model and show that, after the application of a suitable feasibility-based bound tightening procedure, the convex relaxation shares the same optimal value and solution of the non-convex problem.
We also establish that the feasible region of the non-convex problem is a lattice and, through that, a necessary and sufficient condition for the non-emptiness of the feasible region.
\end{abstract}
\section{Introduction}
\label{sec:intro}
In this paper we address the problem of speed planning for a road vehicle on an assigned path in such a way that two distinct and conflicting objectives, namely travel time and energy consumption,
are minimized.
\newline\newline\noindent
We consider a simplified model of an half-car.
We denote with:
\begin{itemize}
\item $M$ the mass and $v$ the speed of the vehicle;
\item $\alpha(s)$ the slope angle of the road, which is a function of the arc-length position;
\item $F_a(t)$ the aerodynamic drag force;
\item $T_f$ and $T_r$ the traction (if positive) or braking (if negative) forces of the front and rear wheels; 
\item $R_f$ and $R_r$ the rolling resistances of the front and rear wheels.

\end{itemize}
We consider negligible the component of the car acceleration orthogonal to the road is negligible.
Then,  the following dynamic equations, obtained from the balance of forces in the direction parallel and orthogonal to the road, can be derived:
\begin{equation}
\label{eqn_dyn_1a}
  \begin{array}{ll}
    M \dot v(t)= T_f(t)+T_r(t)-R_f(t)-R_r(t)-F_a(t) - Mg \sin \alpha(s(t))\\ [8pt]
0=W_f(t)+W_r(t)-Mg \cos \alpha(s(t)).
\end{array}    
\end{equation}
From the second equation the overall load is $W_f(t)+W_r(t)=Mg \cos \alpha(s(t))$.
The traction (or braking) forces must satisfy constraints
\[
|T_f| \leq \mu W_f, |T_r| \leq \mu W_r,
\]
in order to avoid slipping (we assume that the front and rear axes have the same friction coefficient $\mu$).
From the second equation in~(\ref{eqn_dyn_1a}), we have that
\begin{equation}
  \label{eqn_no_slip}
|T_f(t)+T_r(t)| \leq \mu M g \cos \alpha(s(t)).
\end{equation}
We assume that the car has four-wheel drive, and that drive and braking forces are optimally balanced.
In this case, condition~(\ref{eqn_no_slip}) is sufficient to prevent slipping on front and rear wheels.
We assume that the slope angle $\alpha$ is not too large, and we approximate $\cos \alpha(s(t)) \simeq 1$ in condition~(\ref{eqn_no_slip}).
Setting $F(t)= T_f(t)+ T_r(t)$ as the overall traction and braking force, and $F_r(t)=R_f(t)+R_r(t)$ as the overall rolling resistance, we rewrite the first of~\eqref{eqn_dyn} as
\begin{equation}
  \label{eqn_dyn}
M \dot v(t)= -F_a(t) + F(t) - M g \sin \alpha(s(t)) - F_r.
\end{equation}
The aerodynamic drag is $F_a(t)=\frac{1}{2} \rho A_f c_d v(t)^2= \Gamma  v(t)^2$, where $\rho$ is the air density, $A_f$ is the car cross-sectional area, and $c_d$ is the dimensionless drag coefficient.
We model the overall tires rolling resistance with the simplified model
$F_r(t)= (W_f(t)+W_r(t)) c= M g \cos \alpha(s) c \simeq M g c$, where $c$ is the dimensionless rolling resistance coefficient.
Note that $F_r(t)=F_r$ is independent of $t$.
The power consumed (if positive) or recovered (if negative) by the vehicle is given by 
\[
P(t)= \max\{\eta F(t),F(t)\} v(t).
\]
According to the value of $\eta \in [0, 1]$, we obtain different power consumption behaviors: in case of a thermal engine, the vehicle is not recovering energy when braking ($\eta = 0$).
Otherwise, in case of hybrid or electric engines ($\eta > 0$), the regenerative braking system of the vehicle enables energy recovery during deceleration, which accounts for a negative value for the consumed power, meaning that power is recovered.
\newline\newline\noindent
If we consider as variable the kinetic energy as a function of position $s$, that is $w(s)=\frac{1}{2} v^2(s)$, essentially it is like considering the squared velocity (see, for instance,~\cite{verscheure09}) scaled by a constant factor.
Then, in this way, $w'(s)=v(s) \dot v(s) v(s)^{-1}=\dot v(s)$, and~\eqref{eqn_dyn} becomes
\[
M w'(s)= -F_a(s) + F(s) - M g \sin \alpha(s) - F_r,
\]
where the forces are now considered as functions of distance $s$.
We can then write our optimization problem as follows
\begin{align}
& \min_w \int_{s=0}^{s_f} \left(\lambda \max\{\eta F(s),F(s)\} + \frac{1}{\sqrt{w(s)}}\right)ds		& \nonumber \\
\textrm{such that} &											& \nonumber \\
& w(s) \leq \wmax(s)											& \text{for } s \in [0, s_f], \nonumber \\
& |F(s)| \leq M g \mu											& \text{for } s \in [0, s_f], \nonumber \\
& F(s) \leq \frac{\Pmax}{\sqrt{w(s)}}								& \text{for } s \in [0, s_f], \nonumber \\
& M w'(s)= -F_a(s) + F(s) - M g \sin \alpha(s)	 - F_r					& \text{for } s \in [0, s_f], \nonumber \\
& w(s) \geq 0												& \text{for } s \in [0, s_f], \nonumber \\
& w(0) = \winit,\ w(s_{f})=\wfin.												& \nonumber
\end{align}
Note that:
\begin{itemize}
\item an upper bound $\wmax(s)$ is imposed for the speed $w(s)$, depending on the position along the path, since, e.g., the maximum allowed speed is different along a curve and along a straight road;
\item an upper bound $\mu M g$ is imposed for the traction or braking force;
\item an upper bound $\Pmax$ is imposed for the power;
\item the initial and final speed are fixed to $\winit$ and to $\wfin$, respectively (boundary conditions);
\item the objective function is a weighted sum of the two terms:  
$$
\int_{s=0}^{s_f} \max\{\eta F(s),F(s)\}ds\ \mbox{(energy consumption),} \ \ \ \int_{s=0}^{s_f} \frac{1}{\sqrt{w(s)}}ds \ \mbox{(travel time)}.
$$
The parameter $\lambda\geq 0$ is the weight assigned to the energy consumption term.
It represents the time cost (in seconds) per Joule (its dimension is \unit{\second\per\joule}).
\end{itemize}
We discretize the problem by observing variable and input functions at positions $\{0,h,2h,\ldots,(n-1)h\}$,  where $h$ is the discretization step and $n \in \mathbb{N}$ is such that $(n-1)h \approx s_f$ is the length of the path.
We also approximate derivatives with finite differences and the integral in the objective function
with the Riemann sum of the intervals.
Then, the resulting discretized version of the previous problem is:
\begin{align}
& \min_w \sum_{i=1}^{n-1} \left(\lambda \max\{\eta F_i,F_i\} + \frac{1}{\sqrt{w_i}}\right)h		& \nonumber \\
\textrm{such that} &										& \nonumber \\
& w_i \leq \wmax_i										& \text{for } i \in \{1, \ldots, n\}, \nonumber \\
& |F_i| \leq M g \mu										& \text{for } i \in \{1, \ldots, n-1\}, \nonumber \\
& F_i \leq \frac{\Pmax}{\sqrt{w_i}}							& \text{for } i \in \{1, \ldots, n-1\}, \nonumber \\
& \frac{M}{h} (w_{i+1}-w_i) = -\Gamma w_i+ F_i - M g(\sin \alpha_i - c)	& \text{for } i \in \{1, \ldots, n-1\}, \nonumber \\
& w_i \geq 0											& \text{for } i \in \{1, \ldots, n\}, \nonumber \\
& w_1 = \winit,\ w_n=\wfin.											& \nonumber
\end{align}
After setting $f_i=\frac{F_i}{M}$, we can also write the problem as follows ($\gamma=\frac{\Gamma}{M}$):
\begin{align}
& \min_w \sum_{i=1}^{n-1} \left(\lambda M \max\{\eta f_i,f_i\} + \frac{1}{\sqrt{w_i}}\right)h		& \nonumber \\
\textrm{such that} &										& \nonumber \\
& w_i \leq \wmax_i										& \text{for } i \in \{1, \ldots, n\}, \nonumber \\
& |f_i| \leq g \mu										& \text{for } i \in \{1, \ldots, n-1\}, \nonumber \\
& f_i \leq \frac{\Pmax}{M\sqrt{w_i}}							& \text{for } i \in \{1, \ldots, n-1\}, \nonumber \\
& \frac{1}{h} (w_{i+1}-w_i) = -\gamma w_i+ f_i - g(\sin \alpha_i - c)	& \text{for } i \in \{1, \ldots, n-1\}, \nonumber \\
& w_i \geq 0											& \text{for } i \in \{1, \ldots, n\}, \nonumber \\
& w_1 = \winit,\ w_n=\wfin.											& \nonumber
\end{align}
\subsection{Literature Review} 
On the other hand, if $\lambda$ is large, this problem becomes a minimum-fuel one.
We review the literature both on minimum-time and minimum-fuel problems, since they are strongly related.
When $\lambda=0$, the problem simplifies to a minimum-time optimization.
Conversely, for large values of $\lambda$, it turns into a minimum-fuel problem.
Given the close relationship between these two formulations, we examine the literature on both.

\subsubsection{Minimum-time speed planning}
For instance, these constraints can be velocities, accelerations and maximal steering
In some cases, the second choice considerably simplifies the resulting optimization problem.
Reference~\cite{RaiPerCGL17jerk} solves the problem with an additional jerk constraint through a heuristic approach that computes a speed profile by bisection; the proposed method is very efficient, however, the optimality of the obtained profile is not guaranteed.
Several studies have investigated the control problem for minimizing travel time along a predefined route.
Some approaches define vehicle speed as a function of time, while others express it in terms of arc-length position along the path.
The latter approach often leads to a more straightforward optimization process.
Among studies adopting the time-dependent representation,~\cite{MunOllPraSim:94,munoz1998speed} introduce an iterative method based on third-degree polynomial concatenation.
Reference~\cite{SolNun2006} develops a solution using the five-spline scheme from~\cite{bianco2006optimal}, whereas~\cite{Villagra-et-al2012} proposes a closed-form speed profile constructed from three distinct classes of speed profiles.
Additionally,~\cite{CheHeBuHanZha2014} presents an algorithm that generates a piecewise linear speed profile.
For studies that define speed as a function of arc-length,~\cite{Minari16,minSCL17} explore a problem that, after discretization using finite elements, can be reformulated into a convex optimization model.
This allows for efficient computational approaches, as demonstrated in~\cite{CLMNV19,minSCL17,LippBoyd2014}.
Notably,~\cite{minSCL17} describes an algorithm with linear-time complexity concerning the number of variables, providing an optimal solution through spatial discretization.
The method divides the path into $n$ equal segments, approximating the speed derivative via finite differences.
Some works tackle this problem directly in continuous time~\cite{consolini2020solution,Frego16,Velenis2008}.
In particular,~\cite{consolini2020solution} derives an exact continuous-time solution without relying on finite-dimensional reductions, offering a highly efficient method suitable for real-time speed planning.
Reference~\cite{RaiPerCGL17jerk} incorporates a jerk constraint and employs a heuristic bisection approach to compute a speed profile.
While this method is computationally efficient, it does not guarantee optimality.

\subsubsection{Fuel consumption minimization}
Other studies focus on determining an optimal speed profile that minimizes fuel consumption over a given route.
Many works examine energy-efficient speed planning for metro trains.
For example,~\cite{wang2011optimal} formulates the optimal energy-saving speed profile as a Mixed Integer Linear Programming (MILP) problem.
In~\cite{kang2011ga}, Kang et al. introduce an algorithm that optimizes train speed profiles by controlling coasting points, employing Genetic Algorithms (GA).
Reference~\cite{calderaro2014algorithm} identifies an optimal sequence of fundamental control regimes to minimize energy consumption, considering track topology (slopes and curves), vehicle mechanics, power system characteristics, and the influence of regenerative braking.
\newline\newline\noindent
For private vehicles, reducing fuel consumption poses a significant challenge, as drivers typically lack awareness of the most energy-efficient speed trajectories.
The globally optimal speed profile depends on multiple factors and requires substantial computational effort to determine~\cite{gustafsson2009automotive,russell2002integrated}.
However, recent advancements in communication technology, sensor integration, and in-vehicle computing have made real-time optimization more practical.
Some public transport systems in Europe, for instance, already exchange information with traffic lights~\cite{koenders2008cooperative}, while in the USA, researchers are exploring the broadcast of red-light timings to improve safety~\cite{intersections2008cooperative}.
Several algorithms have been proposed to enhance speed optimization, such as predicting optimal speed profiles when approaching traffic signals~\cite{asadi2010predictive} or minimizing energy consumption over routes that include traffic lights~\cite{ozatay2012analytical}.
However, many existing methods rely on expensive onboard computational resources, limiting their real-time feasibility.
Cloud computing has emerged as a promising alternative for handling real-time speed optimization~\cite{wollaeger2012cloud}.
Building on this concept,~\cite{ozatay2014cloud} introduces a real-time Speed Advisory System (SAS) that leverages cloud computing to generate optimal speed profiles based on traffic and geographic data, thus enabling a global optimization framework.
A more aggressive energy-saving strategy is the ``pulse-and-glide'' (PnG) technique~\cite{li2015effect,li2016fuel,li2015mechanism,xu2015fuel}, which consists of rapid acceleration to a target speed followed by a coasting phase to a lower speed.
PnG has been shown to achieve substantial fuel savings in vehicles equipped with continuously variable transmissions (CVT)~\cite{li2015mechanism} and step-gear mechanical transmissions~\cite{xu2015fuel}, with theoretical analyses confirming its optimality and efficiency.
Reference~\cite{kim2019real} proposes a real-time PnG algorithm in the speed-acceleration domain, achieving fuel savings of 3\%--5\% while maintaining near-optimal performance.
Other studies focus on the role of road slope in fuel efficiency~\cite{hellstrom2009look,kamal2011ecological,musardo2005ecms,xu2018design}.
The optimization of speed and control strategies to account for terrain variations significantly influences overall energy consumption.
Model Predictive Control (MPC) methods, such as those proposed by Kamal et al.~\cite{kamal2011ecological} and Hellstr\"om et al.~\cite{hellstrom2009look}, enhance fuel efficiency but impose high computational costs due to system nonlinearities.
An alternative approach, the Equivalent Consumption Minimization Strategy (ECMS)~\cite{musardo2005ecms}, reduces computational complexity while delivering near-optimal fuel savings by relying on instantaneous slope data.
Both predictive and non-predictive strategies are viable for connected automated vehicles (CAVs), yet achieving a balance between fuel efficiency and computational feasibility remains a key challenge for practical deployment.

\subsection{Outline of the paper}
In Section~\ref{sec:convrel} we introduce a convex relaxation of the problem.
In Section~\ref{sec:exact} we briefly sketch the exactness result for this relaxation (under suitable assumptions) when we drop
the boundary condition $w_n=\wfin$ at step $n$, already proved in our paper \cite{ARDIZZONI2025}.
Moreover, we show that the relaxation is not exact any more as soon as we include also the boundary condition.
In Section~\ref{sec:feasboundtighten} we introduce a feasibility-based bound tightening technique, while in Section~\ref{sec:exactfin} we prove that the addition of the lower limit for the variables obtained through the bound tightening procedure makes the convex relaxation exact.
Next, in Section~\ref{sec:furtherbound} we introduce further feasibility-based bound tightening techniques.
In Section~\ref{sec:lattice} we prove that the feasible region of our problem is a lattice and from that we derive a necessary and sufficient condition to establish its non-emptiness, also discussing a procedure to verify such condition, based on the iterated application of the bound tightening techniques discussed in Sections~\ref{sec:feasboundtighten} and~\ref{sec:furtherbound}.
Finally, in Section~\ref{sec:concl} we draw some conclusions.

\section{Convex relaxation of the problem}
\label{sec:convrel} 
The constraints of our problem are the following:
\begin{align}
& \frac{M}{\Pmax\sqrt{w_{i}}}\geq  \frac{1}{h} (w_{i+1}-w_i) + \gamma w_i + g (\sin \alpha_i + c) & i \in \{1, \ldots, n-1\} \label{eq:powermax} \\[8pt]
& \frac{1}{h} (w_{i+1}-w_i) + \gamma w_i + g (\sin \alpha_i + c) \leq g\mu &  i \in \{1, \ldots, n-1\} \label{eq:forcemax1} \\[8pt]
& \frac{1}{h} (w_{i+1}-w_i) + \gamma w_i + g (\sin \alpha_i + c)\geq -g\mu & i \in \{1, \ldots, n-1\} \label{eq:forcemax2} \\[8pt]
& w_i \leq \wmax_i & i \in \{1, \ldots, n\} \label{eq:speedmax} \\[8pt]
& w_i \geq 0 &  i \in \{1, \ldots, n\} \label{eq:nonneg} \\[8pt]
& w_1 = \winit, \ w_n=\wfin. & \label{eq:limitcond}
\end{align}
We rewrite constraints~(\ref{eq:speedmax})--(\ref{eq:limitcond}) as follows:
\begin{align}
& w_i \leq \wmax_i & i \in \{1, \ldots, n\} \label{eq:wmax} \\[8pt]
& w_i \geq \wmin_i & i \in \{1, \ldots, n\} \label{eq:wmin}, 
\end{align}
where $\wmin_i=0$, for $i\in \{2,\ldots,n-1\}$, while $\wmin_1=\wmax_1=\winit$ and $\wmin_n=\wmax_n=\wfin$ (note that the two boundary conditions~(\ref{eq:limitcond}) are split into the constraints $\winit\leq w_1\leq \winit$ and $\wfin\leq w_n\leq \wfin$).
Moreover, after introducing the variables $t_i,f_i$, $i=1,\ldots,n-1$, we can replace the maximum power constraints~(\ref{eq:powermax}) with the following constraints:
\begin{align}
& t_i = \frac{1}{\sqrt{w_{i}}}&  i \in \{1, \ldots, n-1\} \label{eq:pow1} \\[8pt] 
& t_i \geq \frac{M f_i}{\Pmax} &  i \in \{1, \ldots, n-1\} \label{eq:pow2} \\[8pt]  
& f_i = \frac{1}{h} (w_{i+1}-w_i) + \gamma w_i + g (\sin \alpha_i + c) &  i \in \{1, \ldots, n-1\}. \label{eq:pow3}
\end{align}
We denote with $X$ the set of points $w\in \mathbb{R}^n$, $t,f\in \mathbb{R}^{n-1}$ fulfilling the constraints~(\ref{eq:speedmax})--(\ref{eq:limitcond}) with the maximum power constraints~(\ref{eq:powermax}) replaced by the constraints
(\ref{eq:pow1})--(\ref{eq:pow3}).
Such region is non-convex due to the equality constraints $t_i=\frac{1}{\sqrt{w_{i}}}$.
If we replace such equality constraints with the (convex) inequality constraints  $t_i\geq \frac{1}{\sqrt{w_{i}}}$,
we end up with a convex region $X'\supset X$.
Then, our problem is the following (non-convex) problem
\begin{equation}
\label{eq:probfix}
\min_{(w,f,t)\in X} \sum_{i=1}^{n-1} h \left(\lambda M\max\{\eta f_i,f_i\} + t_i\right),
\end{equation}
and its convex relaxation
\begin{equation}
\label{eq:relaxfix}
\min_{(w,f,t)\in X'} \sum_{i=1}^{n-1} h \left(\lambda M\max\{\eta f_i,f_i\} + t_i\right).
\end{equation}
Note that such convex relaxation can be formulated as a Second Order Cone Programming (SOCP) problem.
Indeed, after introducing the variables $z_i,y_i$, $i=1,\ldots,n-1$, we can rewrite the constraints $t_i \geq w_{i+1}^{-\frac{1}{2}}$, as
\[
1 \leq z_i y_i,\quad y_i^2 \leq t_i,\quad z_i^2 \leq t_i w_{i},
\]
(see, e.g., ~\cite{alizadehSecondorderConeProgramming2003}).
These are SOCP constraints.
In general, any constraint of form $a^2 \leq b c$, with $a,b,c$ positive is equivalent to
\[
\left\| \begin{array}{cc} 2 a \\ b-c\end{array} \right\| \leq b+c.
\]

\section{Exactness result for the case without final boundary condition}
\label{sec:exact} 
First of all, we define the critical (squared) speed $\bar w$ as the solution of
\begin{equation}\label{eq:critical_speed}
\frac{\Pmax}{M\sqrt{\bar w}} = g\mu.
\end{equation}
Next, we introduce a couple of assumptions.
\begin{assumption}
\label{assum:1}
It holds that:
\begin{equation}
\label{ineq:h_suff_cond}
(1 - h\gamma) \bar w - hg (1 + c) > \left(\frac{\Pmax h}{2M(\lambda \gamma \Pmax h + 1 - \lambda)}\right)^{\frac{2}{3}}.
\end{equation}
\end{assumption}
Note that Assumption~\ref{assum:1} always holds if the discretization step $h$ is chosen small enough.
\begin{assumption}
\label{assum:2}
It holds that:
\begin{equation}
\label{eq:condbarw}
\min_{i\in I} \left[\frac{\Pmax}{M}\sqrt{\frac{1 - h\gamma}{\bar w + hg(\sin \alpha_i + c)}}-\frac{\gamma}{1-h\gamma}(\bar w + hg(\sin \alpha_i + c))- g (\sin \alpha_i + c)\right]\geq 0,
\end{equation}
where
$$
I=\{i\in \{1,\ldots,n\}\ :\ \wmax_i>\bar{w}\}.
$$
\end{assumption}
Now, let us consider~\eqref{eq:probfix} and its convex relaxation~(\ref{eq:relaxfix}) where in both problems we remove the boundary condition
$w_n=\wfin$.
It is possible to prove that in this case, under Assumptions~\ref{assum:1} and~\ref{assum:2}, the two problems share the same optimal value and solution.
\begin{proposition}
\label{prop:condwbar}
Let us assume that Assumptions~\ref{assum:1} and~\ref{assum:2} hold.
Then, the optimal solution of the relaxed problem~(\ref{eq:relaxfix}) is feasible and, thus, optimal for problem~\eqref{eq:probfix} if we remove constraint $w_n=\wfin$ in both problems.
\end{proposition}
We refer to \cite{ARDIZZONI2025} for a detailed proof of Proposition~\ref{prop:condwbar} but here we briefly sketch the main ideas of the proof.
If we denote with $(w^\star,t^\star,f^\star)$ the optimal solution of the convex relaxation~(\ref{eq:relaxfix}) without the boundary condition $w_n=\wfin$, then this is feasible and, thus, optimal for problem~(\ref{eq:probfix}),
again without the boundary condition $w_n=\wfin$, if and only if
$$
t_i^\star=\frac{1}{\sqrt{w_{i}^\star}}, \ \ \ \forall\ i\in \{1,\ldots,n-1\}.
$$ 
The proof of this fact can be split into two parts:
\begin{itemize}
\item show that there does not exist $i\in \{1,\ldots,n-2\}$ such that
\begin{equation}
\label{eq:contr1}
t_i^\star>\frac{1}{\sqrt{w_{i}^\star}}\ \ \ \mbox{and}\ \ \ t_{i+1}^\star=\frac{1}{\sqrt{w_{i+1}^\star}},
\end{equation}
i.e., it cannot happen that the maximum power constraint is violated by the optimal solution of the relaxation at step $i$, and is not violated at step $i+1$;
\item show that 
\begin{equation}
\label{eq:contr2}
t_{n-1}^\star>\frac{1}{\sqrt{w_{n-1}^\star}},
\end{equation}
cannot hold, i.e., the last maximum power constraint cannot be violated.
\end{itemize}
Both results are proved by contradiction.
Assuming that~(\ref{eq:contr1}) holds, we prove that one can always build a new feasible solution of~(\ref{eq:relaxfix}) with lower objective function value
by simply reducing by a small amount $\delta$ the value of $w_{i+1}^\star$, thus contradicting
optimality of $(w^\star,t^\star,f^\star)$.
In a completely analogous way we can prove that~(\ref{eq:contr2}) cannot hold.
In this case the objective function value can be reduced after reducing by a small amount $\delta$ the value of $w_{n}^\star$.
\newline\newline\noindent
The question now is whether we can extend the result of Proposition~\ref{prop:condwbar} to prove that the optimal solution of problem~(\ref{eq:relaxfix}) is always feasible and, thus, optimal for problem~(\ref{eq:probfix}), also when we include the boundary condition $w_n=\wfin$.
Unfortunately, the answer is no.
However, we can give a characterization of the cases when the result does not hold and, through such characterization, we can discuss how to modify
problem~(\ref{eq:relaxfix}) to make it an exact convex relaxation of~(\ref{eq:probfix}).
By the same argumentation as in the proof of Proposition~\ref{prop:condwbar}, we can show that~(\ref{eq:contr1}) cannot hold.
However,~(\ref{eq:contr2}) may hold.
We cannot apply the proof by contradiction as before since in this problem we cannot decrease the value of $w_n^\star$, which is fixed.
Thus, if the optimal solution of problem~(\ref{eq:relaxfix}) is unfeasible for~(\ref{eq:probfix}), then there exists some $r\geq 1$ such that
\begin{equation}
\label{eq:condviol}
t_i^\star=\frac{1}{\sqrt{w_{i}^\star}}, \ \forall i\in \{1,\ldots,n-r-1\}\ \ \ \mbox{and}\ \ \ t_{i}^\star>\frac{1}{\sqrt{w_{i}^\star}}\ \forall i\in \{n-r,\ldots,n-1\}.
\end{equation}
Stated in another way, if the optimal value of~(\ref{eq:relaxfix}) differs (is lower than) the optimal value of problem~(\ref{eq:probfix}), then the optimal solution of problem~(\ref{eq:relaxfix}) violates all the last $r$ maximum power constraints for some $r$.

\section{Feasibility-based bound tightening}
\label{sec:feasboundtighten}
First, let us introduce the function:
$$
\ell(w)=(1-h\gamma)w+h\min\left\{\frac{\Pmax}{M\sqrt{w}},g\mu \right\}=\left\{
\begin{array}{ll}
(1-h\gamma)w+hg\mu & w\leq \bar{w} \\ [8pt]
(1-h\gamma)w+h\frac{\Pmax}{M\sqrt{w}} & w> \bar{w} .
\end{array}
\right.
$$
Next, we introduce the following assumption.
\begin{assumption}
\label{assum:3}
It holds that:
\begin{equation}
\label{ineq:incrass}
1 - h\gamma -h\frac{\Pmax}{2M\bar{w}^{\frac{3}{2}}} \geq 0.
\end{equation}
\end{assumption}
Note that this assumption is fulfilled for $h$ small enough.
Under this assumption we can prove the following lemma.
\begin{lmm}
\label{lem:incr}
Under Assumption~\ref{assum:3}, function $\ell$ is increasing for $w>0$.
\end{lmm}
\begin{proof}
Under Assumption~\ref{assum:3} it holds that $1-h\gamma>0$ so that $\ell$ is increasing for $w\leq \bar{w}$.
For $w>\bar{w}$ we have that the derivative of $\ell$ is:
$$
\ell'(w)=1-h\gamma -h\frac{\Pmax}{2Mw^{\frac{3}{2}}}.
$$
Since $\ell'$ is increasing with respect to $w$, if $\ell'(\bar{w})\geq 0$, then $\ell'(w)>0$ for all $w>\bar{w}$, i.e., $\ell$ is increasing for $w>\bar{w}$.
Now it is enough to observe that
$$
\ell'(\bar{w})=1 - h\gamma -h\frac{\Pmax}{2M\bar{w}^{\frac{3}{2}}},
$$
so that, under Assumption~\ref{assum:3}, $\ell$ is increasing also for $w>\bar{w}$.
\end{proof}
Now, let us consider
the region $X_1\supset X$ defined as follows:
\begin{equation}
\label{eq:powforcons}
X_1=\left\{w\ :\  \frac{1}{h} (w_{i+1}-w_i) + \gamma w_i + g (\sin \alpha_i + c)\leq \min\left\{\frac{\Pmax}{M\sqrt{w_{i}}},g\mu \right\},\ i=1,\ldots,n-1,\ w\geq \wmin \right\},
\end{equation}
i.e., the region defined by the maximum power and maximum force constraints~(\ref{eq:powermax}) and~(\ref{eq:forcemax1}), and by the lower limit constraints~(\ref{eq:wmin}).
\newline\newline\noindent
For a fixed value $\bar{w}_{i+1}$, $i\in \{1,\ldots,n-1\}$, let us consider the $i$-th maximum power and maximum force constraints and the $i$-th nonnegative constraint:
\begin{equation}
\label{eq:powforcons1}
\begin{array}{l}
 \frac{1}{h} (\bar{w}_{i+1}-w_i) + \gamma w_i + g (\sin \alpha_i + c)\leq\frac{\Pmax}{M\sqrt{w_i}} \\ [8pt]
 \frac{1}{h} (\bar{w}_{i+1}-w_i) + \gamma w_i + g (\sin \alpha_i + c)\leq g \mu \\ [8pt]
 \wmin_i \leq w_i.
 \end{array}
 \end{equation}
 This system of inequalities can be rewritten as:
 \begin{equation}
\label{eq:powforcons2}
\begin{array}{l}
\bar{w}_{i+1} + hg (\sin \alpha_i + c)\leq \ell(w_i)  \\ [8pt]
 \wmin_i \leq w_i.
 \end{array}
 \end{equation}
We denote with $\xi_1^{\wmin}(\bar{w}_{i+1})$ the smallest value which can be assigned to $w_i$ so that
the $i$-th maximum power, maximum force and nonnegative constraints are fulfilled, i.e, so that~(\ref{eq:powforcons2}) is fulfilled.
We prove the following lemma.
\begin{lmm}
\label{lem:xi}
We have that $\xi_1^{\wmin}(\bar{w}_{i+1})=\wmin_i$ if $\ell(\wmin_i)\geq \bar{w}_{i+1} + hg (\sin \alpha_i + c)$, while for  $\ell(\wmin_i)< \bar{w}_{i+1} + hg (\sin \alpha_i + c)$, we have that
$\xi_1^{\wmin}(\bar{w}_{i+1})>\wmin_i$ and $\xi_1^{\wmin}$ is increasing with $\bar{w}_{i+1}$.
\end{lmm} 
\begin{proof}
In view of Lemma~\ref{lem:incr}, $\ell$ is increasing for $w>0$.
Then, if $\ell(\wmin_i)\geq \bar{w}_{i+1} + hg (\sin \alpha_i + c)$, we have that  $\xi_1^{\wmin}(\bar{w}_{i+1})=\wmin_i$, while
for  $\ell(\wmin_i)< \bar{w}_{i+1} + hg (\sin \alpha_i + c)$, there exists a unique value $\xi_1^{\wmin}(\bar{w}_{i+1})>\wmin_i$ which can be assigned to $w_i$ such that
$\ell(\xi_1^{\wmin}(\bar{w}_{i+1}))= \bar{w}_{i+1} + hg (\sin \alpha_i + c)$.
Moreover, the value
$\xi_1^{\wmin}(\bar{w}_{i+1})$ is increasing with $\bar{w}_{i+1}$.
\end{proof}
Now we define a point $l$ as follows:
\begin{equation}
\label{eq:recurtildew}
l_n=\wfin,\ \ l_{j}=\xi_1^{\wmin}(l_{j+1})\ \ \ j\in \{1,\ldots,n-1\}.
\end{equation}
The following proposition shows that each value $l_{j}$ is the smallest value which can be attained by variable $l_{j}$ at points $w\in X_1$.
\begin{proposition}
\label{prop:smallest}
For $j=1,\ldots,n$, the smallest possible value of variable $w_{j}$ at points belonging to $X_1$ is $l_{j}$.
\end{proposition}
\begin{proof}
We prove this by induction.
The result is obviously true for $j=n$ since in $X_1$ we have that $w_n$ is bounded from below by $\wfin$.
Now , for $j< n$, let us assume that the result is true for all $i\in \{j+1,\ldots,n\}$ and let us prove it for $j$.
In particular, the inductive assumption guarantees that
$l_{j+1}$ is the smallest value of variable $w_{j+1}$ which can be assigned to such variable at points in $X_1$.
Moreover, for each $w\in X_1$ for which $w_{j+1}=\bar{y}$, where $\bar{y}$ is a fixed value, the smallest value which can be attained by $w_{j}$ is $\xi_1^{\wmin}(\bar{y})$, i.e., the smallest solution
of the system of inequalities~(\ref{eq:powforcons1}) or the equivalent system~(\ref{eq:powforcons2})  for $\bar{w}_{j+1}=\bar{y}$.
In view of Lemma~\ref{lem:xi}, the minimum value for $\xi_1^{\wmin}(\bar{y})$ is attained when $\bar{y}$ is smallest, i.e., when $\bar{y}=l_{j+1}$, since, by the inductive assumption $l_{j+1}$ is the smallest possible value of $w_{j+1}$ at points in $X_1$.
Thus, the smallest possible value for $w_{j}$ at points in $X_1$ is $\xi_1^{\wmin}(l_{j+1})=l_{j}$, as we wanted to prove.
\end{proof}
Since $X_1 \supset X$, we can also conclude that the values $l_j$, $j=1,\ldots,n$, are valid lower bounds for the variables at all points in $X$, i.e., all feasible solutions of problem~(\ref{eq:probfix}).
In other words, we can say that the recursive relation~(\ref{eq:recurtildew}) defines a {\em feasibility-based bound tightening} procedure, i.e., a procedure that employs some of the constraints defining the region $X$, namely 
the maximum power and maximum force constraints~(\ref{eq:powermax}) and~(\ref{eq:forcemax1}) together with the lower limit constraints~(\ref{eq:wmin}, to tighten the (lower) bounds on the variables.
\newline\newline\noindent
The following proposition introduces a condition under which problem~(\ref{eq:probfix}) does not admit any feasible solution.
\begin{proposition}
\label{prop:unfeas}
If for some $i\in \{1,\ldots,n\}$ it holds that $l_i>\wmax_i$, then $X=\emptyset$.
\end{proposition}
\begin{proof}
In view of Proposition~\ref{prop:smallest},  the smallest possible value of $w_{i}$ at points in $X_1$ is $l_{i}$.
As observed before, since
$X_1\supset X$, then also at points in $X$ the value of $w_{i}$ cannot be lower than $l_{i}$.
But since 
$l_{i}>\wmax_{i}$, while at points in $X$ we must have that $w_{i}\leq \wmax_{i}$, we can conclude that the feasible region of~(\ref{eq:probfix}) is empty.
\end{proof}
Now, let us assume that $l\leq \wmax$.
In the next section we prove that the relaxation becomes exact once we add the new lower bounds for the variables.

\section{Exact convex relaxation with final boundary condition}
\label{sec:exactfin}
According to Proposition~\ref{prop:smallest},
inequalities
\begin{equation}
\label{eq:validlim}
w_i\geq l_i, \ \ \ i=1,\ldots,n,
\end{equation} 
can be added to~(\ref{eq:probfix}) without modifying its feasible region and, consequently, its optimal value.
We can also add them to the convex relaxation~(\ref{eq:relaxfix}) but in this case
they might modify the feasible region and the optimal value.
More precisely, we can prove the following proposition, stating that after the addition of~(\ref{eq:validlim}), the optimal value of~(\ref{eq:relaxfix}) becomes equal to the optimal value of~(\ref{eq:probfix}).
\begin{proposition}
\label{prop:exactrel}
Let $(w^\star,t^\star,f^\star)$ be the optimal solution of problem~(\ref{eq:relaxfix}) with the additional constraints~(\ref{eq:validlim}).
Then, its optimal value is equal to the optimal value
of~(\ref{eq:probfix}) and the two problems share the same optimal solution.
\end{proposition}
\begin{proof}
To prove this fact, we can show that~(\ref{eq:contr1}) and~(\ref{eq:contr2}) cannot hold.
\newline\newline\noindent
It is easily seen that~(\ref{eq:contr2}) cannot hold.
Indeed, we have that $w_n^\star=\wfin$.
Then, by definition of $l_{n-1}$, the
$(n-1)$-th maximum power and force constraints~(\ref{eq:powermax}) and~(\ref{eq:forcemax1}) are fulfilled by any $w$ such that $w_{n-1}\geq l_{n-1}$ and, thus, they are fulfilled by $w^\star$.
Consequently, we have that $t_{n-1}^\star=\frac{1}{\sqrt{w_{n-1}^\star}}$, thus contradicting~(\ref{eq:contr2}).
\newline\newline\noindent
To prove that also
~(\ref{eq:contr1}) cannot hold, we can proceed by contradiction as in the proof of Proposition
~\ref{prop:condwbar}: (i) first, we assume, by contradiction, that~(\ref{eq:contr1}) holds for some $i\in \{1,\ldots,n-2\}$; (ii) then, we
 reduce by a small amount $\delta>0$ the value $w_{i+1}^\star$, (iii) finally, we show that we are led to a contradiction, since the objective function value of the new solution is lower than the objective function value
 of $(w^\star,t^\star,f^\star)$, which contradicts optimality.
 However, some care is needed, since now we are also imposing lower limits for variables and the reduction of $w_{i+1}$ may not be possible if $w_{i+1}^\star=l_{i+1}$, i.e., if $w_{i+1}^\star$ is equal to the lower limit of $w_{i+1}$.
 Actually, what we can prove is that $w_{i+1}^\star=l_{i+1}$ is not possible if~(\ref{eq:contr1}) holds.
Indeed, in case $w_{i+1}^\star=l_{i+1}$, any value $w_i\geq l_i$ would be such that
 the $i$-th  maximum power and force constraints~(\ref{eq:powermax}) and~(\ref{eq:forcemax1}) would be fulfilled, so that $t_i^\star=\frac{1}{\sqrt{w_{i}^\star}}$, which contradicts~(\ref{eq:contr1}).
\end{proof}

\section{Further bound-tightening techniques}
\label{sec:furtherbound}
Although the feasibility-based bound tightening procedure described in Section~\ref{sec:feasboundtighten} is enough for the definition of new lower bounds for the problem variables which make 
the convex relaxation exact for problem~(\ref{eq:probfix}), we can further tighten bounds
of the variables as follows.
\newline\newline\noindent
Let us introduce the set $X_2\supset X$ defined as follows: 
$$
X_2=\left\{w\ :\  \frac{1}{h} (w_{i+1}-w_i) + \gamma w_i + g (\sin \alpha_i + c)\leq \min\left\{\frac{\Pmax}{M\sqrt{w_{i}}},g\mu \right\},\ \ i=1,\ldots,n-1,\ w\leq \wmax\right\},
$$
i.e., the region defined by the maximum power and maximum force constraints~(\ref{eq:powermax}) and~(\ref{eq:forcemax1}), 
and by the upper limit constraints~(\ref{eq:wmax}).
\newline\newline\noindent
For a fixed value $\bar{w}_{i}$, $i\in \{1,\ldots,n-1\}$, let us consider the $i$-th maximum power and maximum force constraints and the $(i+1)$-th maximum speed constraint:
$$
\begin{array}{l}
 \frac{1}{h} (w_{i+1}-\bar{w}_{i}) + \gamma \bar{w}_{i} + g (\sin \alpha_i + c)\leq\frac{\Pmax}{M\sqrt{\bar{w}_i}} \\ [8pt]
 \frac{1}{h} (w_{i+1}-\bar{w}_{i}) + \gamma \bar{w}_{i} + g (\sin \alpha_i + c)\leq g \mu \\ [8pt]
 w_{i+1} \leq \wmax_{i+1}.
 \end{array}
 $$
 Then, a valid upper bound for variable $w_{i+1}$ is:
 $$
 w_{i+1}\leq \xi_2^{\wmax}(\bar{w}_i)=\min\left\{\wmax_{i+1},\ell(\bar{w}_i)-h g (\sin \alpha_i + c)\right\}.
 $$
Similarly to the proof of Proposition~\ref{prop:smallest}, we can show that the recursive relation
$$
u_1=\winit,\ \ u_{j+1}=\xi_2^{\wmax}(u_j)\ \ \ j\in \{1,\ldots,n-1\},
$$
defines a set of valid upper bounds for the variables $w_j$.
\newline\newline\noindent
Next, let us introduce the set $X_3\supset X$ defined as follows: 
$$
X_3=\left\{w\ :\ -\frac{1}{h} (w_{i+1}-w_i) - \gamma w_i - g (\sin \alpha_i + c)\leq g\mu,\ \ i=1,\ldots,n-1,\ w\geq \wmin \right\},
$$
i.e., the region defined by the minimum force constraints~(\ref{eq:forcemax2}),
and by the lower limit constraints~(\ref{eq:wmin}).
\newline\newline\noindent
For a fixed value $\bar{w}_{i}$, $i\in \{1,\ldots,n-1\}$, let us consider the $i$-th minimum force constraint and the $(i+1)$-th lower limit constraint:
$$
\begin{array}{l}
 -\frac{1}{h} (w_{i+1}-\bar{w}_{i}) - \gamma \bar{w}_{i} - g (\sin \alpha_i + c)\leq g \mu \\ [8pt]
 \wmin_{i+1}\leq w_{i+1} .
 \end{array}
$$
 Then, a valid lower bound for variable $w_{i+1}$ is:
 $$
 w_{i+1}\geq \xi_3^{\wmin}(\bar{w}_i)=\max\left\{\wmin_{i+1},(1-h\gamma)\bar{w}_i-h g (\sin \alpha_i + c+\mu)\right\}.
 $$
Similarly to the proof of Proposition~\ref{prop:smallest}, we can show that the recursive relation
$$
l'_1=\winit,\ \ l'_{j+1}=\xi_3^{\wmin}(l_j)\ \ \ j\in \{1,\ldots,n-1\},
$$
defines a set of valid lower bounds for the variables $w_j$.
\newline\newline\noindent
Finally, let us introduce the set $X_4\supset X$ defined as follows: 
$$
X_4=\left\{ w\ :\ -\frac{1}{h} (w_{i+1}-w_i) - \gamma w_i - g (\sin \alpha_i + c)\leq g\mu,\ \ i=1,\ldots,n-1,\ w\leq \wmax\right\},
$$
i.e., the region defined by the minimum force constraints~(\ref{eq:forcemax2}), 
and by the upper limit constraints~(\ref{eq:wmax}).
For a fixed value $\bar{w}_{i+1}$, $i\in \{1,\ldots,n-1\}$, let us consider the $i$-th minimum force constraint and the $i$-th maximum speed constraint:
$$
\begin{array}{l}
 -\frac{1}{h} (\bar{w}_{i+1}-w_{i}) - \gamma w_{i} - g (\sin \alpha_i + c)\leq g \mu \\ [8pt]
 w_i\leq \wmax_{i} .
 \end{array}
$$
Then, a valid upper bound for variable $w_{i}$ is:
$$
w_{i}\leq \xi_4^{\wmax}(\bar{w}_{i+1})=\min\left\{\wmax_i,\frac{\bar{w}_{i+1}+h g (\sin \alpha_i + c+\mu)}{1-h\gamma}\right\}.
$$
Similarly to the proof of Proposition~\ref{prop:smallest}, we can show that the recursive relation
$$
u'_n=\wfin,\ \ u'_{j}=\xi_4^{\wmax}(u'_{j+1})\ \ \ j\in \{1,\ldots,n-1\},
$$
defines a set of valid upper bounds for the variables $w_j$.
\newline\newline\noindent
Note that $X=X_1\cap X_2\cap X_3\cap X_4$.

\section{Lattice structure and a necessary and sufficient condition for feasibility}
\label{sec:lattice}
We denote with $\vee$ and $\wedge$ the component-wise minimum and maximum of two vectors of dimension $n$, respectively, i.e., given $w,w'\in \mathbb{R}^n$, we have that for each $i=1,\ldots,n$:
$$
[w\vee w']_i=\min\{w_i,w'_i\},\ \ \ [w\wedge w']_i=\max\{w_i,w'_i\}.
$$
\begin{proposition}
\label{prop:lattice}
Let Assumption~\ref{assum:3} hold.
If $w, w'\in X$, then $w\wedge w'\in X$ and $w\vee w'\in X$, i.e., $(X,\wedge,\vee)$ is a lattice.
\end{proposition}
\begin{proof}
We notice that $w, w'\in X$ implies that $w\wedge w'\in X$ and $w\vee w'\in X$ fulfill the lower limit constraints~(\ref{eq:wmin}), the upper limit constraints~(\ref{eq:wmax}), and the boundary conditions
(\ref{eq:limitcond}).
Therefore, we only need to prove that $w\wedge w'\in X$ and $w\vee w'\in X$ also fulfill the maximum power constraints~(\ref{eq:powermax}), and the maximum and minimum force constraints~(\ref{eq:forcemax1})--(\ref{eq:forcemax2}).
Also recalling the definition of $\ell$, we can rewrite these constraints as follows:
\begin{equation}
\label{eq:reform}
\begin{array}{ll}
w_{i+1}\leq \ell(w_i)-h g (\sin \alpha_i + c) & i=1,\ldots,n-1 \\ [8pt]
w_{i+1}\geq (1-h\gamma)w_i +h g (\sin \alpha_i + c+\mu)  & i=1,\ldots,n-1,
\end{array}
\end{equation}
where the first inequality is equivalent to the $i$-th maximum power constraint~(\ref{eq:powermax}) and maximum force constraint~(\ref{eq:forcemax1}).
Let us assume, w.l.o.g., that $w_{i+1}\vee w'_{i+1}=w_{i+1}$ and $w_{i+1}\wedge w'_{i+1}=w'_{i+1}$.
\newline\newline\noindent
If $w_i\vee w'_i=w_i$ and $w_i\wedge w'_i=w'_i$, then constraints~(\ref{eq:reform}) are obviously fulfilled by 
$w\vee w'$ and by $w\wedge w'$.
\newline\newline\noindent
Therefore, let us see what happens when $w_i\vee w'_i=w'_i$ and $w_i\wedge w'_i=w_i$.
First of all, we observe that:
$$
\begin{array}{l}
w_{i+1}\vee w'_{i+1}=w_{i+1}\leq w'_{i+1}\leq \ell(\overbrace{w'_i}^{w_{i}\vee w'_{i}})-h g (\sin\alpha_i + c) \\ [8pt]
w_{i+1}\vee w'_{i+1}=w_{i+1}\geq (1-h\gamma)w_i +h g (\sin\alpha_i + c+\mu)\geq (1-h\gamma) \underbrace{w'_{i}}_{w_{i}\vee w'_{i}} +h g (\sin\alpha_i + c+\mu),
\end{array}
$$
that is, constraints~(\ref{eq:reform}) are fulfilled by $w\vee w'$.
Next, recalling the increasing monotonicity of $\ell$ established in Lemma~\ref{lem:incr} under Assumption~\ref{assum:3}, we observe that:
$$
\begin{array}{l}
w_{i+1}\wedge w'_{i+1}=w'_{i+1}\leq \ell(w'_i)-h g (\sin\alpha_i + c)\leq \ell(\overbrace{w_i}^{w_{i}\wedge w'_{i}})-h g (\sin\alpha_i + c) \\ [8pt]
w_{i+1}\wedge w'_{i+1}=w'_{i+1}\geq w_{i+1}\geq (1-h\gamma)\underbrace{w_i}_{w_{i}\wedge w'_{i}} +h g (\sin\alpha_i + c+\mu),
\end{array}
$$
i.e., constraints~(\ref{eq:reform}) are fulfilled by $w\wedge w'$.
\end{proof}
Since $(X,\vee,\wedge)$ is a lattice, $X\neq \emptyset$ implies that the component-wise maximum and minimum over $X$, i.e., the points $z',y'$ defined as follows:
 $$
z'_i=\max_{w\in X} w_i,\ \ \ y'_i=\min_{w\in X} w_i, \ i=1,\ldots,n,
$$
are such that $z',y'\in X$.
But how can we establish whether $X\neq \emptyset$ and, in case it is, how can we compute $z',y'$?
We proceed as follows.
Let
$$
z_i=\max_{w\in X_2\cap X_4} w_i,\ \ \ y_i=\min_{w\in X_1\cap X_3} w_i,\ i=1,\ldots,n,
$$
i.e.:
\begin{itemize}
\item $z_i$ is the maximum value which can be attained by $w_i$ if we ignore the lower bound constraints on the variables, i.e., the non-negativity constraints for variables $w_i$, $i\in \{2,\ldots,n-1\}$ together with $w_1\geq \winit$
and $w_n\geq \wfin$;
\item $y_i$ is the minimum value which can be attained by $w_i$ if we ignore the maximum speed constraints (where $\wmax_1=\winit$ and $\wmax_n=\wfin$).
\end{itemize}
Note that $X_1\cap X_3$ and $X_2\cap X_4$ are both nonempty lattices, so that $z\in X_1\cap X_3$ and $y\in X_2\cap X_4$ (later on we will introduce a procedure for their computation).
We prove the following proposition stating a necessary and sufficient condition for $X\neq \emptyset$.
\begin{proposition}
\label{prop:empty}
We have that $X\neq \emptyset$ if and only if $y\leq z$.
 Moreover, if $X\neq \emptyset$, then $z'=z$ and $y'=y$.
\end{proposition}
\begin{proof}
If $y\leq z$, then $\wmax_i\geq z_i\geq y_i\geq 0$, for $i\in \{2,\ldots,n-1\}$, and:
$$
\begin{array}{lll}
\wmax_1=\winit\geq z_1\geq y_1\geq \winit & \Rightarrow & z_1=y_1=\winit \\
\wmax_n=\wfin\geq z_n\geq y_n\geq \wfin & \Rightarrow & z_n=y_n=\wfin,
\end{array}
$$
so that $y,z\in X$ and $X\neq \emptyset$.
\newline\newline\noindent
Moreover, since $X\subset X_1\cap X_3, X_2\cap X_4$, in this case $z$ and $y$ are the component-wise maximum and minimum over $X$, respectively, i.e., $z'=z$ and $y'=y$.
\newline\newline\noindent
Instead, if $y\not\leq z$, we can conclude that $X=\emptyset$.
Indeed, since each $z_i$ is an upper limit for $w_i$ over $X_2\cap X_4$ and, thus, over $X$, while each $y_i$ is a lower limit for $w_i$ over $X_1\cap X_3$ and, thus, over $X$,
if for some $i$ it holds that $y_i>z_i$, we can conclude that $X=\emptyset$.
\end{proof}
The next question is how to compute $z,y$.
We proceed as follows.
First, we introduce the functions:
$$
\begin{array}{ll}
B_1(l)=p &\mbox{where } p_n=l_n,\ p_j=\xi_1^l(p_{j+1}),\ \ j=1,\ldots,n-1 \\ [8pt]
B_2(u)=p &\mbox{where }  p_1=u_1,\ p_{j+1}=\xi_2^l(p_{j}),\ \ j=1,\ldots,n-1 \\ [8pt]
B_3(l)=p &\mbox{where }  p_1=l_1,\ p_{j+1}=\xi_3^l(p_{j}),\ \ j=1,\ldots,n-1 \\ [8pt]
B_4(u)=p &\mbox{where }  p_n=u_n,\ p_j=\xi_4^u(p_{j+1}),\ \ j=1,\ldots,n-1.
\end{array}
$$
Then, we can apply Algorithm~\ref{alg:buildsol} to compute $z,y$.
More precisely, we show that Algorithm~\ref{alg:buildsol} is an iterative algorithm converging to $z,y$ or establishing that $X=\emptyset$.
The algorithm works as follows.
At line~\ref{lin:1} the vectors $u^1$ and $l^1$ are initialized with the vectors $\wmax$ and $\wmin$, respectively.
Next, at iteration $k$ we compute: 
\begin{itemize}
\item at line~\ref{lin:2} a new upper bound vector $u^{k+\frac{1}{2}}$ for the variables $w$.
By definition of $B_2$, it holds that $u^{k+\frac{1}{2}}\leq u^k$;
\item at line~\ref{lin:3} a new upper bound vector $u^{k+1}$ for the variables $w$.
By definition of $B_4$, it holds that $u^{k+1}\leq u^{k+\frac{1}{2}}$;
\item at line~\ref{lin:4} a new lower bound vector $l^{k+\frac{1}{2}}$ for the variables $w$.
By definition of $B_1$, it holds that $l^{k+\frac{1}{2}}\geq l^k$;
\item at line~\ref{lin:5} a new lower bound vector $l^{k+1}$ for the variables $w$.
By definition of $B_3$, it holds that $l^{k+1}\geq l^{k+\frac{1}{2}}$.
\end{itemize}
Then, we check whether the distance between $u_k$ and $u_{k+1}$, and the one between $l^k$ and $l^{k+1}$ is within a given small threshold distance $\varepsilon>0$.
If yes, we stop
and return the vectors $u^{k+1}$ and $l^{k+1}$.
Otherwise, if $u^{k+1}\not\geq l^{k+1}$, since $z\leq u^{k+1}$ and $l^{k+1}\leq y$, we have that $z\not\geq y$ and, in view of Proposition~\ref{prop:empty}, we can stop
establishing that $X=\emptyset$.
\newline\newline\noindent
The following proposition proves a convergence result for Algorithm~\ref{alg:buildsol}.
\begin{proposition}
If we set $\varepsilon=0$, then Algorithm~\ref{alg:buildsol} either stops after a finite number of iteration returning $z,y$, or generates two infinite sequences $\{u^k\}$ and $\{l^k\}$, the former converging from above to $z$, the latter converging from below to $y$.
\end{proposition}
\begin{proof}
First we observe that the two sequences $\{u^k\}$ and $\{l^k\}$ are nonincreasing and nondecreasing, respectively.
Moreover, for each $k$ it holds that $u^k\geq u^{k+\frac{1}{2}}\geq u^{k+1}\geq z$ and  $l^k\leq l^{k+\frac{1}{2}}\leq l^{k+1}\leq y$.
By definition of $u^{k+\frac{1}{2}}, u^{k+1}, l^{k+\frac{1}{2}}, l^{k+1}$, we have that for $i\in \{1,\ldots,n-1\}$:
\begin{equation}
\label{eq:convcons}
\begin{array}{l}
\frac{1}{h} (u^{k+\frac{1}{2}}_{i+1}-u^{k+\frac{1}{2}}_{i}) + \gamma u^{k+\frac{1}{2}}_{i} + g (\sin \alpha_i + c)\leq\frac{\Pmax}{M\sqrt{u^{k+\frac{1}{2}}_{i}}} \\ [8pt]
\frac{1}{h} (u^{k+\frac{1}{2}}_{i+1}-u^{k+\frac{1}{2}}_{i}) + \gamma u^{k+\frac{1}{2}}_{i} + g (\sin \alpha_i + c)\leq g \mu \\ [8pt]
-\frac{1}{h} (u^{k+1}_{i+1}-u^{k+1}_{i}) - \gamma u^{k+1}_{i}- g (\sin \alpha_i + c)\leq g \mu \\ [8pt]
\frac{1}{h} (l^{k+\frac{1}{2}}_{i+1}-l^{k+\frac{1}{2}}_{i}) + \gamma l^{k+\frac{1}{2}}_{i} + g (\sin \alpha_i + c)\leq\frac{\Pmax}{M\sqrt{l^{k+\frac{1}{2}}_{i}}} \\ [8pt]
\frac{1}{h} (l^{k+\frac{1}{2}}_{i+1}-l^{k+\frac{1}{2}}_{i}) + \gamma l^{k+\frac{1}{2}}_{i} + g (\sin \alpha_i + c)\leq g \mu \\ [8pt]
-\frac{1}{h} (l^{k+1}_{i+1}-l^{k+1}_{i}) - \gamma l^{k+1}_{i}- g (\sin \alpha_i + c)\leq g \mu.
 \end{array}
\end{equation}
If at iteration $k$ we have that $\|u^{k+1}-u^k\|=\|l^{k+1}-l^k\|=0$, then $u^k=u^{k+\frac{1}{2}}=u^{k+1}$ and  $l^k=l^{k+\frac{1}{2}}=l^{k+1}$, and, in view of~(\ref{eq:convcons}), $u^k\in X_2\cap X_4$, while $l^k\in X_1\cap X_3$, so that $z=u^k$ and $y=l^k$.
\newline\newline\noindent
Instead, if $\|u^{k+1}-u^k\|=\|l^{k+1}-l^k\|=0$ never occurs, then we observe that the sequences $\{u^k\}$ and $\{l^k\}$ are monotonic and, if we never stop, they are also bounded, since $u_k$ cannot fall below $\wmin$, while $l^k$ cannot be larger than $\wmax$.
Therefore, the sequences converge to $\bar{u}\geq z$ and $\bar{l}\leq y$, respectively.
Taking the limit for $k\rightarrow \infty$ in~(\ref{eq:convcons}), we have that $\bar{u}\in X_2\cap X_4$, while $\bar{l}\in X_1\cap X_3$, so that $z=\bar{u}$ and $y=\bar{l}$.
\end{proof}
\begin{algorithm}[H]
\small
\begin{algorithmic}[1] 
\Procedure{ComputeZY}{$\wmax,\wmin,\varepsilon>0$}
\State Set $u^1=\wmax$, $l^1=\wmin$, $k=1$, $stop={\tt false}$ \label{lin:1}\;
\While{$stop={\tt false}$}
\State Set $u^{k+\frac{1}{2}}=B_2(u^k)$ \label{lin:2}\;
\State Set $u^{k+1}=B_4(u^{k+\frac{1}{2}})$ \label{lin:3}\;
\State Set $l^{k+\frac{1}{2}}=B_1(l^k)$\label{lin:4}\;
\State Set $l^{k+1}=B_3(l^{k+\frac{1}{2}})$\label{lin:5}\;
\If{$\|u^{k+1}-u^k\|\leq \varepsilon\ \ \mbox{{\tt and}}\ \ \|l^{k+1}-l^k\|\leq \varepsilon$ \label{lin:6}}
\State ${\it stop}={\tt true}$\;
\State \Return{$u^{k+1}, l^{k+1}$}\;
\ElsIf{$u^{k+1}\not\geq l^{k+1}$\label{lin:7}}
\State{$X=\emptyset$}
\Else
\State Set $k=k+1$\;
\EndIf
\EndWhile
\EndProcedure
\caption{\label{alg:buildsol} Procedure to compute the vectors $z,y$ or to establish that $X=\emptyset$}
\end{algorithmic}
\end{algorithm}

\section{Conclusions}
\label{sec:concl}
In this paper we addressed the problem of planning the speed of a vehicle along an assigned path by minimizing a weighted sum of energy consumption and travel time.
The resulting mathematical model is a non-convex problem.
The paper extends the result already proved in our paper \cite{ARDIZZONI2025}, where we proved that the problem without a boundary condition about the final speed of the vehicle can be solved through a convex relaxation of its mathematical model.
However, the result cannot be extended to the case when the final speed of the vehicle is fixed.
We introduced a way to strengthen the convex relaxation through a feasibility-based bound tightening technique which allows to impose lower limits for the problem variables, and we proved that the strengthened convex relaxation is exact (i.e., it shares the same optimal value and solution of the non-convex problem).
We also studied the properties of the feasible region of the non-convex problem, proving that it is a lattice.
Through such properties we have also been able to derive a necessary and sufficient condition to establish whether the feasible region is empty or not.

\bibliographystyle{plain} 
\bibliography{biblio}
\end{document}